\documentclass{article}

\usepackage{amsmath,amssymb,a4wide,amsthm}

\usepackage{graphicx}

\newtheorem{theorem}{Theorem}
\newtheorem{lemma}[theorem]{Lemma}
\newtheorem{proposition}[theorem]{Proposition}

\newtheorem{definition}[theorem]{Definition}
\newtheorem{remark}[theorem]{Remark}

\numberwithin{theorem}{section}
\numberwithin{equation}{section}

\title{Area limit laws for symmetry classes of staircase polygons}

\author{\sc Christoph Richard\dag, Uwe Schwerdtfeger\dag and Bhalchandra 
Thatte\ddag\\
\normalsize
\dag Fakult\"at f\"ur Mathematik, Universit\"at Bielefeld,\\
\normalsize
Postfach 10 01 31, 33501 Bielefeld, Germany\\
\normalsize
\ddag Mathematics and Computer Science Building, University of Canterbury,\\ 
\normalsize
Private Bag 40800, Christchurch, New Zealand}

\begin{document}

\maketitle

\begin{abstract}
We derive area limit laws for the various symmetry classes 
of staircase polygons on the square lattice, in a uniform 
ensemble where, for fixed perimeter, each polygon occurs 
with the same probability. This complements a previous study 
by Leroux and Rassart, where explicit expressions for the 
area and perimeter generating functions of these classes 
have been derived.
\end{abstract}

\smallskip


\section{Introduction}

Models of planar lattice polygons have a rich history. The
most challenging member is self-avoiding polygons \cite{MS93}, 
a model of interest not only in enumerative combinatorics, 
but also in the natural sciences such as physics and chemistry. 
Some solvable subclasses, obtained by imposing  convexity or 
directedness constraints, have been understood in detail, see 
\cite{BM96} for an overview. An important example on the square 
lattice are \textit{staircase polygons}, also called parallelogram 
polygons. These have been enumerated with respect to various 
parameters such as perimeter, area and generalisations thereof, 
site-perimeter, and radius of gyration, see \cite{DGV87,BM96,R06,L07} 
and references therein. For some of these parameters, their 
distribution has been asymptotically analysed. An example is the 
area distribution in a uniform ensemble where, for fixed perimeter, 
all polygons occur with the same probability. We will call this 
ensemble the \textit{uniform fixed perimeter ensemble}. For a large 
class of planar polygon models including staircase polygons, the 
\emph{Airy distribution} emerges as the limit law of area \cite{D99,R05}. 
There is compelling numerical evidence that the Airy distribution 
also appears as the limit distribution of area in self-avoiding 
polygons \cite{RGJ01,R07}.

When counting polygons, translated copies of a polygon are identified, 
but rotated or reflected versions of a polygon are often distinguished. 
Motivated by applications such as benzenoid counting in chemistry 
\cite{VG02}, one may be led to identify objects that are related by 
point symmetries. The task arises to enumerate classes of polygons up 
to symmetries of the underlying lattice. Two different counting problems 
occur. One may ask for the number of polygons which are fixed by a given 
subgroup of the lattice symmetries. One may also count the number of 
polygon orbits with respect to a such a subgroup. These two counting 
problems are related. In fact, it is possible to express the generating 
function for orbit counts in terms of those for fixed point counts, via 
the Lemma of Burnside. This has been studied for various classes of 
column-convex polygons by Leroux and co-workers \cite{LRR98,LR01,GL05}, 
and explicit expressions for generating functions have been obtained.

Symmetry classes of staircase polygons on the square lattice are analysed 
in \cite{LR01}, by two different approaches. One of them uses a bijection 
between staircase polygons and Dyck paths due to Delest and Viennot 
\cite{DV84}. By this bijection, some of the polygon symmetry classes can 
be identified with corresponding symmetry classes of Dyck paths. The 
toolbox of Dyck path analysis can then be used to solve the corresponding 
polygon counting problems. A second approach uses the 
Temperley--Bousquet-M\'elou method \cite{BM96} to obtain explicit 
expressions for the perimeter and area generating functions. In addition, 
it is shown in \cite{LR01} that the number of polygons fixed by some 
non-trivial symmetry is asymptotically negligible to the total number of 
polygons. This implies that the number of symmetry orbits is 
asymptotically equal to the total number of polygons.

In this article, we derive the area limit laws for the symmetry subclasses
of polygons, within the uniform fixed perimeter ensemble. Wheras the result
for the full class of staircase polygons is known already \cite{D99,R05}, the 
corresponding problem for the various symmetry subclasses has apparently not 
been studied before. We will use an elementary self-contained approach, which 
is based on a simple decomposition of staircase polygons \cite{R06}. This 
yields, for \emph{every} subclass, a $q$-difference equation for its 
perimeter and area generating function. Explicit expressions could be 
obtained from this equation, which might differ from those of \cite{LR01}. 
We do not focus on solving, but  on manipulating the $q$-difference equation 
in order to derive the limit law, by an application of the moment method. 
Whereas this approach has been used previously in different contexts, see 
e.g.~\cite{T91,T95,FPV98,D99,NT03}, we would like to stress that our 
alternative derivation, based on the method of \emph{dominant balance} 
\cite{R05,R07}, simplifies the analysis. As a consequence, corrections to 
the asymptotic behaviour might be mechanically obtained, compare \cite{R02}. 

Our results extend those of \cite{LR01} to a 
complete symmetry analysis, and to a detailed discussion of asymptotics. 
We remark that the area limit laws cannot easily be obtained from the functional
equations in \cite{LR01}, or from the explicit expressions given there. One can 
use the asymptotic description of Dyck paths by Brownian excursions 
 \cite{A91,B72} to obtain part of our results.
Whereas this requires some results from stochastic processes, 
our approach only needs basic probability theory and singularity
analysis of generating functions, as e.g.~described in \cite{FS07}.

The plan of this paper is as follows. We will discuss 
the various square lattice symmetry subgroups, 
and characterise the associated polygon classes by 
certain decompositions. This induces functional
equations for their generating functions, from which
the area limit laws are obtained, by an application of
the moment method and the method of dominant balance. 
We will also indicate how some of our results may be obtained 
within a stochastic approach, and we finally discuss 
extensions of our results.  

\section{Symmetry classes and functional equations}

We explain the models, introduce basic constructions, 
and set the notation, following \cite{R06}.
We will then derive functional equations for the perimeter and
area generating function of the symmetry subclasses of
staircase polygons.

\smallskip

\begin{figure}[htb]
\begin{center}
\includegraphics[width=80mm]{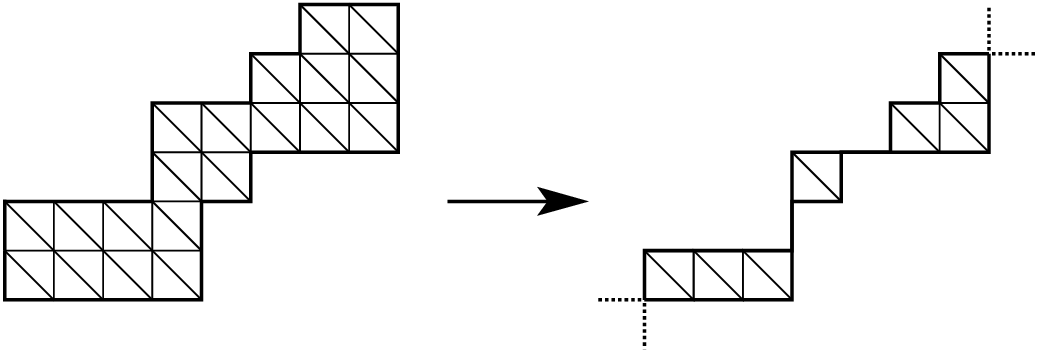}
\end{center}
\caption{\label{fig:rem}
\small The set of staircase polygons is in one-to-one correspondence with 
the set of ordered sequences of edges and staircase polygons. A corresponding
combinatorial bijection is characterised by shifting the upper walk of a 
staircase polygon one unit down and one unit to the right, and by then removing
the first and the last edge of each walk.}
\end{figure}

Consider two fully directed walks on the edges of the square lattice (i.e., 
walks stepping only up or right), which both start at the origin and end 
in the same vertex, but have no other vertex and no edge in common. The edge set of such 
a configuration is called a {\em staircase polygon}, if it is nonempty. 
For a given staircase polygon, 
consider the construction of moving the upper directed walk one unit down and 
one unit to the right. For each walk, remove its first and its last edge.
The resulting object is a sequence of (horizontal and vertical) edges and 
staircase polygons, see Figure \ref{fig:rem}. The unit square yields the
empty sequence. 
It is easy to see that this construction describes a combinatorial bijection
between the set $\mathcal{P}$ of staircase polygons and the set $\cal Q$ of 
ordered sequences of edges and staircase polygons. Let us 
denote the corresponding map by $f:\cal P\to \cal Q$. Thus, for a staircase 
polygon $P\in\cal P$, we have $f(P)=(Q_1,\ldots,Q_n)\in\cal Q$, where $Q_i$ is, 
for $i=1,\ldots, n$, either a single edge or a staircase polygon. 
We denote the single horizontal edge by $e_h$, and the single vertical edge
by $e_v$. The image of the unit square is the empty sequence $n=0$, which we 
occasionally identify with a single point, denoted by $pt.$ A variant of this 
construction will be used below, in order to derive a functional equation for the 
generating function of the staircase polygon symmetry classes.

The perimeter of a staircase polygon $P\in\mathcal{P}$ is defined to be the 
number of its edges. The half-perimeter equals the number of (negative) 
diagonals $n_0(P)$ plus one. The area $n_1(P)$ of a staircase polygon $P$ 
is defined to be the number of its enclosed squares. It equals the sum 
of the lengths of its (negative) diagonals. See Figure \ref{fig:rem} for 
an illustration. The \emph{weight} of a staircase polygon $P$ is the monomial 
$w_P(x,q)=x^{n_0(P)+1}q^{n_1(P)}$. The half-perimeter and area generating 
function of a subclass $\mathcal{C}\subseteq\mathcal{P}$ of staircase 
polygons is the (formal) power series
\[
C(x,q):=\sum_{P\in {\cal C}}w_P(x,q).
\]
Observe that for $e_h$, $e_v$, $pt$, and for $P\in {\cal P}$ we have \cite{R06}
\[
\begin{split}
w_{f^{-1}(pt)}(x,q)=x^2q, \qquad &w_{f^{-1}(e_{h})}(x,q)=x^3q^2, \qquad 
w_{f^{-1}(e_{v})}(x,q)=x^3q^2,
\\
w_{f^{-1}(P)}&(x,q)=x^2q\cdot w_{P}(xq,q).
\end{split}
\]
For a polygon $P\in\mathcal{P}$, consider $f(P)=(Q_1,\ldots,Q_n)$.
In order to retrieve $P$ from $(Q_1,\ldots,Q_n)$, translate $P_{i}=
f^{-1}(Q_{i})$ in such a way that its lower left square coincides 
with the upper right square of $P_{i-1}$, for $i\in\{2,\ldots,n\}$. We say 
that $P$ is the \emph{concatenation} of $(P_1,\ldots,P_n)$, and
write $P=c(P_1,\ldots,P_n)$. The weight $w_P(x,q)$ of $P$ is retrieved 
from the weights of $P_1,\ldots,P_n$ via \cite{R06}
\[
w_{c(P_1,\ldots,P_n)}(x,q)=\frac{1}{(x^2q)^{n-1}}w_{P_1}(x,q)\cdot\ldots\cdot w_{P_n}(x,q).
\]
Denote by $\widetilde{\mathcal{P}} \subseteq \mathcal{P}$ the subset of polygons 
$\widetilde{P}=f^{-1}(P)$, where $P\in {\cal P}\cup\{e_h,e_v\}$. We have
established a combinatorial bijection between the set $\mathcal{P}$ 
and the set of ordered sequences from $\widetilde{\mathcal{P}}$.

\medskip

\begin{figure}
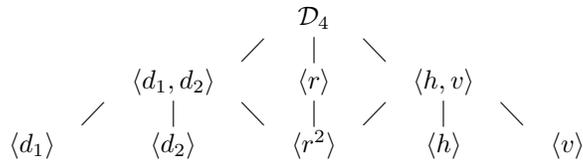
\label{SGL}
\[
\begin{array}{rcccccccl}
& & &  & {\mathcal D}_4&  &  & &  \\
& & & \diagup& \vert & \diagdown & & & \\
& & \langle d_1,d_2\rangle & &\langle r \rangle& & \langle h,v \rangle & & \\
 &\diagup &\vert & \diagdown &\vert &  \diagup& \vert & \diagdown  & \\
\langle d_1\rangle & & \langle d_2\rangle & & \langle r^2 \rangle &  &\langle h\rangle & &\langle v\rangle 
\end{array}
\]

\caption{\small The lattice of subgroups of ${\cal D}_4$. The rotation about $\pi/2$ 
is denoted by $r$, the reflections in the positive and the negative diagonal are 
denoted by $d_1$ and $d_2$, and the reflections in the horizontal and vertical 
axes are denoted by $h$ and $v$. The identity is omitted.} 
\end{figure}

The group of point symmetries of the square lattice is the dihedral group 
$\mathcal{D}_4$. Its non-trivial subgroups are depicted in Figure~\ref{SGL}. 
Note that the above decomposition respects any subgroup of the square
lattice point symmetries. This observation is the key to deriving functional 
equations for the generating functions of the symmetry subclasses.
In the proof of the following proposition, will treat two cases in some detail,
the remaining ones being handled similarly.  

\begin{proposition}\label{prop:feqns}
The half-perimeter and area generating functions of the staircase 
polygon symmetry subclasses satisfy the following functional equations.

\begin{enumerate}
\item Class $\cal P$ of all staircase polygons with generating function $P(x,q)$:
\begin{equation} \label{feqP}
P(x,q)=\frac{x^2q}{ 1- 2xq-P(xq,q)}.
\end{equation}
\item Class $\cal S$ of $\langle r^2\rangle$-symmetric staircase polygons 
with generating function $S(x,q)$:
\begin{equation}\label{feqS}
S(x,q)=\frac{1}{x^2q}\left(1+2xq+S(xq,q)\right)P(x^2,q^2).
\end{equation}
\item Class of $\langle d_1\rangle$-symmetric 
staircase polygons with generating function $D_1(x,q)$:
\[
D_1(x,q)=\frac{x^2q}{1-D_1(xq,q)}.
\]
\item Class of $\langle d_2\rangle$-symmetric 
staircase polygons with generating function $D_2(x,q)$:
\[
D_2(x,q)=\frac{1}{x^2q}\left(1+D_2(xq,q) \right)P(x^2,q^2).
\]
\item Class of $\langle d_1,d_2\rangle$-symmetric staircase polygons with
generating function $D_{1,2}(x,q)$:
\[
D_{1,2}(x,q)=\frac{1}{x^2q}\left(1+D_{1,2}(xq,q) \right)D_1(x^2,q^2).
\]
\item Classes of $\langle h \rangle$-, $\langle v \rangle$-, and 
$\langle h,v \rangle$-symmetric staircase polygons with generating 
function $H(x,q)$:
\[
H(x,q)=x^2qH(xq,q)+x^2q\frac{1+xq}{1-xq}.
\] 
\item Classes $\langle r\rangle$-symmetric staircase 
polygons with generating function $R(x,q)$:
\[
R(x,q)=x^2qR(xq,q)+x^2q.
\]
\end{enumerate}
\end{proposition}
\begin{proof}
Denote the induced group action $\alpha:\mathcal{D}_4\times\mathcal{P}\to\mathcal{P}$
by $\alpha(g,P)=gP$. 

\noindent 
1. The bijection described above implies the following chain of equalities, 
compare also \cite{R06}.
\begin{equation}\label{feqPprf}
\begin{split}
P(x,q) &= \sum_{n=0}^\infty \sum_{(P_1,\ldots,P_n)\in
(\widetilde{\cal P})^n} w_{c(P_1,\ldots,P_n)}(x,q)\\
&= \sum_{n=0}^\infty x^2q \sum_{(P_1,
\ldots,P_n)\in(\widetilde{\cal P})^n} \frac{w_{P_1}(x,q)}
{x^2q}\cdot\ldots\cdot\frac{w_{P_n}(
x,q)}{x^2q}\\
&= x^2q \sum_{n=0}^\infty \left( 
\frac{1}{x^2q}\sum_{P\in\widetilde{\cal P}} 
w_{P}(x,q)\right)^n\\
&= x^2q  \frac{1}{1-\frac{1}{x^2q} \left( w_{f^{-1}(e_h)}(x,q) +w_{f^{-1}(e_v)}(x,q) + 
\sum_{p\in\cal P} w_{f^{-1}(P)}(x,q)\right)
}\\
& =\frac{x^2q}{ 1- 2xq-P(xq,q)}.
\end{split}
\end{equation}

\noindent 
2. For $P\in{\mathcal S}$, we have $f(P)=(P_1,\ldots,P_n,C,r^2P_n,\ldots,
r^2P_1)$, where $C\in{\mathcal S}\cup \{e_v,e_h,pt\}$, and $P_i \in 
\mathcal{P}\cup \{e_v,e_h\}$ for $i=1,\ldots,n$, compare 
Figure~\ref{symstair}.
\begin{figure}\label{symstair}
\begin{center}
\includegraphics[width=80mm]{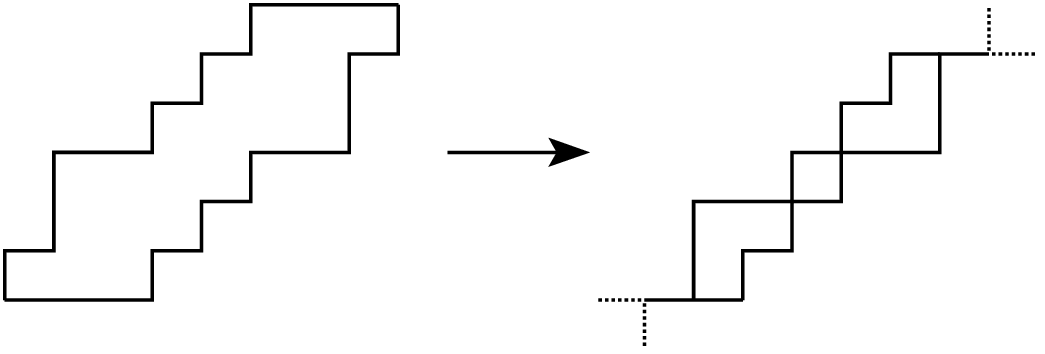}
\caption{$r^2$-symmetric polygon and corresponding sequence of polygons and edges}
\end{center}
\end{figure}
In analogy to the definition of $\widetilde{\mathcal{P}}$ above, define
$\widetilde{\mathcal{S}}\subset\mathcal{P}$ as the pre-image of ${\mathcal S}
\cup \{e_v,e_h,pt\}$ under $f$. Note that concatenation of $Q\in \cal P$ 
with the unit square results in $Q$ again, and that we have $w_Q(x,q)^k=
w_Q(x^k,q^k)$. With $P(x,q)$ as above, this yields
\begin{equation}\label{feqSprf}
\begin{split}
S(x,q) &= \sum_{n=0}^\infty \sum_{(P_1,\ldots,P_n,C)\in
(\widetilde{\cal P})^n\times\widetilde{{\mathcal S}}} w_{c(P_1,\ldots,P_n,C,r^2P_n,\ldots,r^2P_1)}(x,q)\\
&=\sum_{C\in \widetilde{{\cal S}}}w_C(x,q) \sum_{n=0}^\infty  \sum_{(P_1,
\ldots,P_n)\in(\widetilde{\cal P})^n} \frac{w_{P_1}(x,q)^2}
{(x^2q)^2}\cdot\ldots\cdot\frac{w_{P_n}(x,q)^2}{(x^2q)^2}\\
&= \left(w_{f^{-1}(pt)}+w_{f^{-1}(e_h)}+w_{f^{-1}(e_v)}+\sum_{C\in {\cal S}}w_{f^{-1}(C)}\right)
 \sum_{n=0}^\infty \left( 
\sum_{Q\in\widetilde{\cal P}} 
\frac{w_{Q}(x^2,q^2)}{(x^2q)^2}\right)^n\\
&= \left(x^2q+2x^3q^2+x^2qS(xq,q)\right)
 \frac{1}{1-\left(2x^2q^2 +P(x^2q^2,q^2) \right)}\\
&= \frac{1+2xq+S(xq,q)}{x^2q}\cdot
 \frac{x^4q^2}{1-2x^2q^2 -P(x^2q^2,q^2)}\\
&=\frac{1}{x^2q}\left(1+2xq+S(xq,q)\right)P(x^2,q^2)
\end{split}
\end{equation}
where the sum over $Q\in\widetilde{\cal P}$ in the third equation is treated as in 
eqn.~\eqref{feqPprf}. In the last step, we applied eqn.~\eqref{feqP}.

\noindent 
3. For a $\langle d_1\rangle$-symmetric polygon $Q$, we 
have $f(Q)=(P_1,\ldots P_n)$, with a $\langle d_1\rangle$-symmetric $P_i$ 
for $i=1,\ldots,n$. A calculation similar to that in eqn.~\eqref{feqPprf} then 
yields the assertion.

\noindent 
4. For a $\langle d_2\rangle$-symmetric polygon $Q$, we 
have $f(Q)=(P_1,\ldots P_n,C,d_2P_n,\ldots,d_2P_1)$, where 
$P_i \in {\cal P}\cup \{e_h,e_v\}$ for $i=1,\ldots,n$, and where $C$ is a $pt$ 
or $\langle d_2 \rangle$-symmetric. Now a computation similar to that in eqn.~\eqref{feqSprf} 
yields the assertion.

\noindent 
5. For a $\langle d_1,d_2\rangle$-symmetric polygon $Q$, we 
have $f(Q)=(P_1,\ldots P_n,C,d_2P_n,\ldots,d_2P_1)$, with 
$\langle d_1\rangle$-symmetric $P_i$ for $i=1,\ldots,n$, and 
where $C$ is a $pt$ or $\langle d_1d_2\rangle$-symmetric. A computation similar to 
that of eqn.~\eqref{feqSprf} yields the assertion.

\noindent 
6. Staircase polygons are also characterised by the property that 
they contain the lower left and the upper right corner of their smallest bounding 
rectangles. So the only staircase polygons with $\langle h \rangle$- or 
$\langle v \rangle$- symmetry are rectangles.  $f$ maps a rectangle $Q$ either 
to a single rectangle, or to sequences of vertical (horizontal) edges, 
if the width (height) of $Q$ is 1. This results in the above equation.

\noindent  
7. The only admissible polygons are squares. For a given 
half-perimeter, there is exactly one square. If $n>1$, the function $f$ 
maps a square of half-perimeter $2n$ to the square of half-perimeter 
$2n-2$, and it maps the unit square to $pt$. We obtain the claimed 
equation.

\end{proof}

\begin{remark} \rm Equations of the above form appear
in different contexts. Examples are classes
of directed lattice paths, counted by length and area under the
path \cite{NT03}, or classes of simply generated trees \cite{MM78}, 
counted by number of vertices and internal path length. This is due 
to combinatorial bijections between these classes, which we will partly review in 
Section~\ref{sec:stoch}. In the context of polygon models, 
equations appear for Class 1 in \cite{BMV92} and for Class 6
in \cite{R02}, while Class 7 is trivial. Solutions of some equations
may be given in explicit form, compare \cite{BM96,PB95}.
\end{remark}

\section{Area limit laws}

In this section, we derive the limiting area laws for the various symmetry 
subclasses, in the uniform fixed perimeter ensemble. This will be achieved
by an application of the moment method \cite[Sec.~30]{Bill95}. Such an
approach has been used previously \cite{T91,T95,NT03,D99} in similar contexts,
using some involved computations. We will follow a streamlined version, based 
on the method of dominant balance \cite{R05}, which finally allows to obtain 
the limit distribution by a mechanical calculation. In order to give a 
self-contained description of the method, we will treat the two cases $\cal P$ 
and $\cal S$ in detail, and then indicate the analogous arguments for the 
remaining subclasses.

\subsection{Limit law for $\cal P$}\label{limitP}
A $q$-difference equation for the half-perimeter and area
generating function $P(x,q)$ of all staircase polygons was derived
in Proposition \ref{prop:feqns}. For $q=1$, the resulting quadratic 
equation describes the generating function $P(x,1)$ of staircase polygons, 
counted by half-perimeter. The relevant solution is
\begin{equation}\label{D00}
P(x,1)= \frac{1}{4}-\frac{1}{2}\sqrt{1-4x}+\frac{1}{4}(1-4x) .
\end{equation}  

We are interested in the distribution of area within a uniform ensemble where, for
fixed perimeter $2m$, each polygon has the same probability
of occurrence. We introduce the discrete random variable $X_m$
of area by
\begin{equation}\label{tXm}
\mathbb P_m( X_m=n)=\frac{[x^mq^n]P(x,q)}{[x^m]P(x,1)},
\end{equation}
where $[u^k]f(u)$ denotes the coefficient of order $k$ in the
formal power series $f(u)$. In the following, we will asymptotically 
analyse the moments of $ X_m$.
The answer can be expressed in terms of the \emph{Airy distribution}, 
see \cite{FL01,J07,KMM07} for a discussion of its properties.

\begin{definition}\label{def:Airy}
A random variable $Y$ is {\em Airy distributed} \cite{FL01} if
\begin{equation*}
\frac{\mathbb E[Y^k]}{k!} = \frac{\Gamma(\gamma_0)}{\Gamma(\gamma_k)} 
\frac{\phi_k}{\phi_0},
\end{equation*}
where $\gamma_k=3k/2-1/2$, and where $\Gamma(z)$ is the Gamma function.
The numbers $\phi_k$ satisfy for $k\in\mathbb N$ the quadratic recursion
\begin{equation*}
\gamma_{k-1}\phi_{k-1}+\frac{1}{2}\sum_{l=0}^k \phi_l \phi_{k-l}=0,
\end{equation*}
with initial condition $\phi_0=-1$. 
\end{definition}

\begin{remark} \rm
\textit{i)} In the sequel, we shall make frequent use of \emph{Carleman's 
condition}: A sequence of moments $\{M_m\}_{m\in \mathbb{N}}$ with the 
property $\sum_k(M_{2k})^{-1/(2k)}=\infty$ defines a unique random 
variable $X$ with moments $M_m$, cf. \cite{Fel70}.
\\
\textit{ii)} This implies in particular, that $Y$ is uniquely determined by 
the above moment sequence. Explicit expressions can be given for its 
moments, its moment generating function, and its density. The name relates 
to the asymptotic expansion
\begin{displaymath}
\frac{\rm d}{{\rm d}s}\log \mbox{Ai}(s)\sim \sum_{k\ge0}
(-1)^k \frac{\phi_k}{2^k}s^{-\gamma_k} \qquad (s\to\infty),
\end{displaymath}
where $\mbox{Ai}(x)=\frac{1}{\pi}\int_0^\infty \cos(t^3/3+tx){\rm d}t$
is the Airy function. The Airy distribution appears in a variety of 
contexts. In particular, the random variable $Y/\sqrt{8}$ describes 
the Brownian excursion area. 
\end{remark}

We can now state the following result.

\begin{theorem}\label{theo:Dyck}
For staircase polygons of half-perimeter $m$, the area random variables 
$X_m$ eqn.~\eqref{tXm}, appropriately normalised, converge in distribution,
\begin{displaymath}
\frac{ X_m}{m^{3/2}} \stackrel{d}{\to} \frac{Y}{4} 
\qquad (m\to\infty).
\end{displaymath}
where $Y$ is Airy distributed.
We also have moment convergence.
\end{theorem}

\begin{remark}\rm
The previous theorem is a special case of \cite[Thm.~3.1]{D99} 
and \cite[Thm.~1.5]{R05}. In \cite{D99}, a limit distribution
result is stated for certain algebraic $q$-difference equations, 
together with arguments of a proof using the moment method. In 
\cite{R05}, a general multivariate limit distribution result is proved
for certain $q$-functional equations, using the moment method and the method
of dominant balance. The corresponding argument for staircase 
polygons, using the moment method and the method of dominant balance, 
is sketched in \cite{R06}.
\end{remark}

For pedagogic purposes, we will give a proof of the above result 
using the moment method and the method of dominant balance. In 
particular, we show that the moments of $ X_m$, appropriately 
normalised, converge to those of $Y$. Since $Y$ is uniquely determined by its 
moments by Carleman's condition, this implies that the sequence $\{ 
X_m\}_{m\in\mathbb N}$ converges, after normalisation, 
to $Y$ in distribution and for moments, compare \cite[Thm.~4.5.5]{Chu74}.

We will analyse the asymptotic behaviour of the moments in terms of
the singular behaviour of the associated \emph{factorial moment generating 
functions} which are, for $(k,l)\in \mathbb{N}_0^2$, defined by
\begin{displaymath}
P_{k,l}(x):=
\left.\frac{\partial^{k+l}}{\partial q^k\partial x^l}P(x,q)
\right|_{q=1}.
\end{displaymath}
In particular, $P_{0,0}(x)=P(x,1)$ is the 
generating function of staircase polygons, counted by half-perimeter.
Note that these quantities exist as formal power series and have
the same radius of convergence as $P(x,1)$. This is due to the fact that
the area of a polygon is bounded by the square of its perimeter, 
resulting in $\left[ x^m \right]P(x,q)$ being a \emph{polynomial} in 
$q$. The name results from the identity
\begin{equation}\label{eq:EXm}
\mathbb E_m[(X_m)_k]=\frac{[x^m]P_{k,0}(x)}{[x^m]P_{0,0}(x)},
\end{equation}
where $(a)_k=a\cdot (a-1)\cdot\ldots\cdot (a-k+1)$ is the 
lower factorial. The factorial moment generating functions 
$P_{k,l}(x)$ turn out to be algebraic. Explicit expressions 
may be obtained recursively from the functional equation 
eqn.~\eqref{feqP}, by implicit differentiation 
w.r.t.~$x$ and $q$.

\smallskip

We will study the singular behaviour of the factorial moment
generating functions from their defining functional equation,
and will then infer the asymptotic behaviour of the moments. 
This will be achieved in three steps.
We will first prove the existence of a certain local expansion
for each factorial moment generating function about its singularity, 
by an application of the chain rule (or Fa\`a di Bruno's formula). 
Then we will provide an explicit expression for  the leading term
in the expansion, by an application of the method of dominant balance. 
This will finally be analysed in order to obtain the asymptotic behaviour of
the corresponding moment, by methods from singularity analysis 
of generating functions.

\smallskip

We note that steps one and two are usually performed simultaneously,
the corresponding method being nicknamed \emph{moment pumping} \cite{FPV98}.
Our two-step approach uses an exponent guess, which 
might be obtained from an analysis of the first few factorial moment 
generating functions. It is then shown that the guessed exponent value is an upper 
bound on the true exponent, by an application of Fa\`a di Bruno's formula. 
The corresponding calculation is simpler to perform than the usual asymptotic analysis 
of the functional equation. Once an exponent bound has been established, 
the method of dominant balance can be applied. It yields a 
recursion for the coefficients of the leading singular term in the 
factorial moment generating functions. If the recursion reveals non-zero 
coefficients, this proves that the exponent bound is actually an equality,
thereby verifying the initial guess. 

\medskip

The first step of our method is summarised by the following lemma.
For its statement, recall that a function 
$f(u)$ is \emph{$\Delta$-regular} \cite{FFK04} if it is analytic in the 
\emph{indented disc} $\Delta=\Delta(u_c)=\{u:|u|\le u_c+\eta, |\mbox{arg}
(u-u_c)|\ge\phi\}$ for some real numbers $u_c>0$, $\eta>0$ and $\phi$, 
where $0<\phi<\pi/2$. Note that $u_c\notin\Delta$, where we employ the 
convention $\mbox{arg}(0)=0$. The set of $\Delta$-regular functions is 
closed under addition, multiplication, differentiation, and integration. Moreover, 
if $f(u)\ne0$ in $\Delta$, then $1/f(u)$ exists in $\Delta$ and is $\Delta$-regular.

\begin{lemma}\label{expDkl}
For $(k,l)\in\mathbb N_0^2$, the power series $P_{k,l}(x)$
has radius of convergence $1/4$ and is $\Delta(1/4)$-regular. 
It has a locally convergent expansion about $x=1/4$, as in
eqn.~\eqref{D00} for $(k,l)=(0,0)$, and for $(k,l)\neq(0,0)$ 
of the form
\begin{equation}\label{Dkl}
P_{k,l}(x)=\sum_{r=0}^\infty \frac{d_{k,l,r}}{(1-4x)^{3k/2+l-r/2-1/2}}.
\end{equation}
\end{lemma}

\begin{remark} \rm \textit{i)} The exponent $3k/2$ in 
eqn.~\eqref{Dkl} might be guessed from the asymptotic behaviour of
the mean area $\mathbb E_m[X_m]\sim A m^{3/2}$ of a random polygon. 
The mean area is obtained from $P_{0,0}(x)$ and $P_{1,0}(x)$, 
which might be easily extracted from the $q$-difference equation. 
Note that the coefficients $d_{k,l,0}$ might attain 
zero values at this stage. The recursion eqn.~\eqref{recF} given 
below however implies that all
of them are non-zero.  \\
\textit{ii)} The reasoning in the following 
proof may be used to show that all series $P_{k,l}(x)$ are 
algebraic.\\
\textit{iii)} For $\langle r^2 \rangle$-symmetric polygons, 
our proof below will use properties of the derivatives 
\[
\widetilde{P}_{k,l}(x):=
\left.\frac{\partial^{k+l}}{\partial q^k\partial x^l}
\left(P(x^2,q^2)\right)
\right|_{q=1}.
\]
These functions have all radius of convergence $1/2$, are $\Delta(1/2)$-regular,
and have the same type of expansion as the functions $P_{k,l}(x)$ of
the previous lemma. This may be inferred from the previous lemma by the chain rule or,
more formally, by an application of Fa\`a di Bruno's 
formula \cite{Co96}.
\end{remark}

\begin{proof}
By the argument given above, it is seen that all functions $P_{k,l}(x)$ 
have the same radius of convergence. The statement of the theorem 
is true for $P_{0,0}(x)=P(x,1)$, as follows from the 
explicit expression eqn.~\eqref{D00}. For the general case, we argue by 
induction on $(k,l)$, using the total order $\lhd$ defined by 
\begin{displaymath}
(r,s)\lhd(k,l)\Leftrightarrow r+s<k+l \vee \left(r+s=k+l \wedge r<k \right),
\end{displaymath}
chosen to be compatible with the combinatorics of derivatives. 
Define 
\[
H(x,q):=P(x,q)(1-2xq-P(xq,q))-x^2q,
\]
compare eqn.~\eqref{feqP}. 
Fix $(k,l)\rhd(0,0)$. An application of Leibniz' rule yields
\begin{equation}\label{dfeq}
\begin{split}
\frac{\partial^{k+l}}{\partial q^k \partial x^l} \left(H(x,q)+x^2q\right)= \sum_{(0,0)
\unlhd (r,s)\unlhd (k,l)} & \binom{k}{r} \binom{l}{s}  \frac{\partial^{k+l-r-s}}
{\partial q^{k-r} \partial x^{l-s}} \left (P(x,q) \right) \\
  &\cdot \, \frac{\partial^{r+s}}{\partial q^r \partial x^s}\left(1-2xq-P(xq,q)\right)
\end{split}
\end{equation}
In fact, terms corresponding to indices $(r,s)\unlhd(k,l)$ with $s>l$ or $r>k$ 
are zero. For the second derivative on the r.h.s.~of eqn.~\eqref{dfeq}, note 
that by the chain rule
\[
\frac{\partial^{r}}{\partial q^r}\left( P(xq,q) \right )=\sum_{i=0}^r
\binom{r}{i}q^{r-i}\left(\frac{\partial^r}{\partial x^{r-i}\partial q^i }
P\right)(xq,q).
\]
Taking further derivatives w.r.t.~$x$, we may write
\begin{equation}\label{Dpq}
\begin{split}
\frac{\partial^{r+s}}{\partial q^r \partial x^s}\left(P(xq,q) 
\right)& = q^{r+s}\left(\frac{\partial^{r+s}}{\partial q^r \partial 
x^s}P\right)(qx,q)\\
 & +\sum_{(i,j)\lhd(r,s)} \left(\frac{\partial^{i+j}}{\partial q^i \partial x^j}
P\right)(xq,q)\cdot w_{i,j}(x,q),
\end{split}
\end{equation}
for polynomials $w_{i,j}(x,q)$ in $x$ and $q$, which satisfy $w_{i,j}(x,q)\equiv0$ 
if $i<r$.
By inserting eqn.~\eqref{Dpq} into eqn.~\eqref{dfeq} and setting $q=1$, 
one observes that only the $(0,0)$ and the $(k,l)$ summand in 
eqn.~\eqref{dfeq} contribute terms with $P_{k,l}(x).$ The terms involving 
$P_{k,l}(x)$ sum up to
\[
P_{k,l}(x)\left(1-2x-2P_{0,0}(x)  \right)=\sqrt{1-4x}\,P_{k,l}(x),
\]
where we used eqn.~\eqref{D00}.
Now the claimed $\Delta(1/4)$-regularity of $P_{k,l}(x)$ follows from the
induction hypothesis, by the closure properties of $\Delta$-regular 
functions. For the particular singular expansion eqn.~\eqref{Dkl} note
that, by induction hypothesis, each of the remaining terms in the 
summation in eqn.~\eqref{dfeq} has an expansion eqn.~\eqref{Dkl}.
Hence, the most singular exponent is bounded by
\[
\left(\frac{3}{2}(k-r)+(l-s)-\frac{1}{2}\right)+
\left(\frac{3}{2}r+s-\frac{1}{2}\right)
=\frac{3}{2}k+l-1.
\]
We conclude that the leading singular exponent of $P_{k,l}(x)$ is at most
$3k/2+l-1/2$, which yields the desired bound, and thus the 
remaining assertion of the theorem.

\end{proof}

The second and third step of our method yield
a proof of Theorem \ref{theo:Dyck}.

\begin{proof}[Proof of Theorem \ref{theo:Dyck}]

We apply the method of dominant balance \cite{R05} in order
to obtain the limit distribution of area. Its idea consists in first
replacing the factorial moment generating functions,
which appear in the formal expansion of $P(x,q)$ about $q=1$,
by their singular expansion of Lemma \ref{expDkl}, and then
in studying the equation implied by the $q$-difference equation 
eqn.~\eqref{feqP}. We may thus write
\begin{equation}\label{form:dombalDyck}
P(x,q)=\frac{1}{4}+(1-q)^{1/3}\,F\left(\frac{1-4x}{(1-q)^{2/3}},(1-q)^{1/3}\right),
\end{equation}
where $F(s,\epsilon)=\sum_r F_r(s)\epsilon^r$ is a formal 
power series in $\epsilon$ with coefficients $F_r(s)$ being formal Laurent series 
in $s^{1/2}$. The series
\begin{equation}\label{form:F0dyck}
F_0(s)=F(s,0)=\sum_k \frac{d_{k,0,0}}{k!}\cdot\frac{(-1)^k}{s^{3k/2-1/2}}
\end{equation}
is some generating function for the leading coefficients of $P_{k,0}(x)$. 
The coefficients $f_k:=d_{k,0,0}/k!$, in turn, determine the asymptotic 
form of the factorial moments $\mathbb E[(X_m)_k]$ in eqn.~\eqref{eq:EXm}, as we will 
see below. We will use the $q$-difference equation eqn.~\eqref{feqP} to
derive a defining equation for  $F_0(s)$. This will lead to a simple 
quadratic recursion for the numbers $f_k$. Use the above form of $P(x,q)$ 
in the $q$-difference equation, introduce $4x=1-s\epsilon^2,$  $q=1-\epsilon^3$,
and expand the functional equation to second order in $\epsilon$. This yields
a Riccati equation for the generating function $F_0(s)$,
\begin{equation}\label{diffeqF}
\frac{\mathrm{d}}{\mathrm{d} s}F_0(s)+4F_0(s)^2-s=0.
\end{equation}
On the level of coefficients of $F_0(s)$, we obtain the recursion
 \begin{equation}\label{recF}
 \gamma_{k-1}f_{k-1}+4\sum_{l=0}^kf_lf_{k-l}=0,
 \end{equation}
with initial condition $f_0=d_{0,0,0}=-1/2$. We infer from the definition
of the Airy distribution that $f_k=2^{-2k-1} \phi_k$. In particular, all
coefficients $f_k$ are non-zero. Noting that the functions $P_{k,0}(x)$ are 
$\Delta(1/4)$-regular, we thus get by an application of the transfer lemma 
\cite[Thm.~1]{FO90} for the factorial moments of $X_m$ the asymptotic form
\begin{displaymath}
\begin{split}
\frac{\mathbb E_m[(X_m)_k]}{k!} &= 
\frac{1}{k!}\frac{[x^m]P_{k,0}(x)}{[x^m]P_{0,0}(x)}
\sim \frac{1}{k!}\frac{[x^m]d_{k,0,0}(1-4x)^{-(3k/2-1/2)}}
{[x^m]d_{0,0,0}(1-4x)^{1/2}}\\ 
&\sim  \frac{f_{k}}{f_{0}}
\frac{\Gamma(-1/2)}{\Gamma(3k/2-1/2)} m^{3k/2}
= \frac{\phi_k}{\phi_0}\frac{\Gamma(\gamma_0)}{\Gamma(\gamma_k)}
\left(\frac{m^{3/2}}{4}\right)^k
\qquad (m\to\infty).
\end{split}
\end{displaymath}
The previous estimate also shows that the factorial moment $\mathbb 
E_m[(X_m)_k]$ is asymptotically equal to the ordinary moment $\mathbb 
E_m[(X_m)^k]$. It follows with \cite[Thm.~4.5.5]{Chu74} that the sequence 
of random variables $\left \{4 m^{-3/2}X_m\right\}_{m\in \mathbb{N}}$ 
converges in distribution to $Y$, where $Y$ is Airy distributed. The 
above reasoning also implies moment convergence.

\end{proof}

\noindent 
\begin{remark} \rm
If we expand the functional equation to higher order in the above example, 
we obtain at order $\epsilon^{r+2}$ a linear differential equation for 
the function $F_r(s)$, which is the generating function for the numbers 
$d_{k,0,r}$ in the expansion eqn.~\eqref{Dkl}, compare \cite{R02}. 
So we can mechanically obtain corrections to the asymptotic behaviour of the 
factorial moment generating functions, and hence to the moments 
of the limit distribution. 
\end{remark}

\subsection{Limit law for $\cal S$}\label{limitS}

The above strategy can also be followed in order to study 
the area law for the class of $\langle r^2 \rangle$-symmetric 
staircase polygons. The result can be expressed in terms of the 
distribution of area of the Brownian meander, see \cite[Thms.~2,3]{T95} 
and the review \cite{J07}.

\begin{definition}\label{def:mean}
The random variable $Z$ of area of the Brownian meander
is given by 
\begin{equation*}
\frac{\mathbb E[Z^k]}{k!} = \frac{\Gamma(\alpha_0)}{\Gamma(\alpha_k)} 
\frac{\omega_k}{\omega_0}\frac{1}{2^{k/2}},
\end{equation*}
where $\alpha_k=3k/2+1/2$. The numbers $\omega_k$ satisfy 
for $k\in \mathbb{N}$ the quadratic recursion
\begin{equation}
\alpha_{k-1}\omega_{k-1}+\sum_{l=0}^k \phi_l 2^{-l} \omega_{k-l}=0,
\end{equation}
with initial condition $\omega_0=1$, where the numbers $\phi_k$
appear in the Airy distribution.
\end{definition}

\begin{remark} \rm By Carleman's condition, the random variable $Z$ is 
uniquely determined by its moments, and explicit expressions are 
known for the moment generating function and the distribution function.
\end{remark}

We are particularly interested in the derivatives
\begin{displaymath}
S_{k,l}(x):=
\left.\frac{\partial^{k+l}}{\partial q^k\partial x^l}S(x,q)
\right|_{q=1},
\end{displaymath}
where $(k,l)\in\mathbb N_0^2$. 
As above, these series exist as formal power series and have the same
radius of convergence. We have the following lemma.
\begin{lemma}\label{lem:symdyck}
For $(k,l)\in\mathbb N_0^2$, the power series $S_{k,l}(x)$
has radius of convergence $1/2$ and is $\Delta(1/2)$-regular. 
It has a locally convergent expansion about $x=1/2$ of the form
\begin{displaymath}
S_{k,l}(x)=\sum_{r\ge0} \frac{s_{k,l,r}}
{(1-2x)^{3k/2+l-r/2+1/2}}.
\end{displaymath}
\end{lemma}

\begin{remark} \rm The following proof can be used to show that all series
$S_{k,l}(x)$ are algebraic.
\end{remark}

\begin{proof}
The proof is analogous to that of Lemma \ref{expDkl}.
Elementary estimates show that all series $S_{k,l}(x)$
have the same radius of convergence. Setting 
$q=1$ in eqn.~\eqref{feqS}, solving for $S_{0,0}(x)$ and expanding 
about $x=1/2$ yields the assertion for $(k,l)=(0,0)$.
We argue by induction on $(k,l)$, using the total order $\lhd$.  
Fix $(k,l)\rhd(0,0)$. Differentiating eqn.~\eqref{feqS} with Leibniz' Rule gives
\begin{equation}\label{dfeqS}
\begin{split}
\frac{\partial^{k+l}}{\partial q^k \partial x^l}S(x,q) & = 
\sum_{(0,0)\unlhd (r,s)\unlhd (k,l)} \binom{k}{r} \binom{l}{s}  
\frac{\partial^{r+s}}{\partial q^{r} \partial x^{s}} \left 
(\frac{P(x^2,q^2)}{x^2q} \right) \\
  &\cdot \, \frac{\partial^{k+l-r-s}}{\partial q^{k-r} 
\partial x^{l-s}}\left(1+2xq+S(xq,q)\right).
\end{split}
\end{equation}
We argue as in the proof of Lemma \ref{expDkl} 
that only the $(0,0)$ summand on the right hand side of 
eqn.~\eqref{feqS} contributes $(k,l)$ derivatives of $S$, and
that all other derivatives of $S$ of order $(r,s)$ satisfy $(r,s)\lhd(k,l)$. 
Setting $q=1$ in \eqref{feqS} and collecting 
all terms involving $S_{k,l}(x)$ on the left hand 
side gives
\begin{equation}
\begin{split}
\left(1-\frac{\widetilde{P}_{0,0}(x)}{x^2}\right)S_{k,l}(x)=
\frac{\widetilde{P}_{0,0}(x)}{x^2}(1+2x)+\sum_{(r,s)\unlhd(k,l)} & \left. 
\frac{\partial^{r+s}}{\partial q^{r} \partial x^{s}} \left(\frac{P(x^2,q^2)}{x^2q} \right)\right|_{q=1} \\
 & \cdot\left(h_{r,s}(x)+\sum_{(i,j)\lhd (k,l)}a_{i,j}S_{i,j}(x)\right),
\end{split}
\end{equation}
where the $h_{r,s}(x)$ are (at most linear) polynomials, and the 
$a_{i,j}$ are some real coefficients. Note also that the terms 
\[
\left. \frac{\partial^{r+s}}{\partial q^{r} \partial x^{s}} \left(\frac{P(x^2,q^2)}{x^2q} \right)
\right|_{q=1}=\sum_{i,j}\binom{r}{i}\binom{s}{j}\widetilde{P}_{i,j}(x)\frac{c_{i,j}}{x^{2+r-i}q^{1+s-j}}
\]
are $\Delta(1/2)$-regular, with an expansion about $x=1/2$ having the same 
exponents as in eqn.~\eqref{Dkl}, see the remark following Lemma~\ref{expDkl}.
We thus get $\Delta(1/2)$-regularity of $S_{k,l}(x)$ by induction, and by the
closure properties of $\Delta$-regular functions.
For the particular expansion, note that the right hand side has a locally convergent expansion about 
$1/2$ with most singular exponent $3k/2+l$, as the factor
$\left. \frac{\partial^{r+s}}{\partial q^{r} \partial x^{s}} \left(\frac{P(x^2,q^2)}{x^2q} \right)\right|_{q=1}$ has an expansion 
with most singular exponent $3r/2+s-1/2$, and the inner sum 
has by induction an expansion with most singular exponent 
at most $3(k-r)/2+(l-s)+1/2.$ The first factor on the left 
hand side has a locally convergent expansion about $1/2$ starting with
\[
\left(1-\frac{\widetilde{P}_{0,0}(x)}{x^2}\right)=-2\sqrt{2}\sqrt{1-2x}+\mathcal{O}(1-2x)
\qquad (x\to1/2).
\]
Solving for $S_{k,l}(x)$ yields the desired expansion.
\end{proof}

\begin{theorem}\label{theo:r2}
The random variables $X_m^{(sym)}$ of area of 
$\langle r^2\rangle$-symmetric staircase polygons of half-perimeter $m$, 
appropriately normalised, converge in distribution,
\begin{displaymath}
\frac{X_m^{(sym)}}{m^{3/2}} \stackrel{d}{\to} \frac{Z}{2} \qquad (m\to\infty),
\end{displaymath}
where $Z$ is the meander area random 
variable. We also have moment convergence.
\end{theorem}

\begin{proof}
We apply the method of dominant balance.
According to Lemma \ref{lem:symdyck}, the generating function
$S(x,q)$ may be expressed as
\begin{equation}\label{form:dombalsymDyck}
S(x,q)=\frac{1}{(1-q)^{1/3}}\,
G\left( \frac{1-2x}{(1-q)^{2/3}},(1-q)^{1/3} \right),
\end{equation}
where $G(s,\epsilon)=\sum G_r(s) \epsilon^r$ is a formal 
power series in $\epsilon$ and $s^{-1/2}$, and 
\begin{displaymath}
G(s,0)=G_0(s) = \sum_{k=0}^\infty  \frac{s_{k,0,0}}{k!}\frac{(-1)^k}{s^{3k/2+1/2}}
\end{displaymath}
is a generating function for the leading coefficients in the
singular expansions of the functions $S_{k,0}(x)$. The functional
equation eqn.~(\ref{feqS}) induces a recursion on the numbers 
$g_k:=s_{k,0,0}/k!$, which determines the limit distribution, as we will 
see below. We insert eqn.~(\ref{form:dombalsymDyck}) together 
with eqn.~(\ref{form:dombalDyck}) into the functional equation, 
introduce $q=1-\epsilon^3$ and $2x=1-s\epsilon^2$, 
and expand the functional equation to order zero in $\epsilon.$
This gives the linear inhomogeneous first order differential equation
\begin{equation}\label{diffeqG}
\frac{\mathrm{d}}{\mathrm{d} s} G_0(s)+4\cdot 2^{1/3}
F_0\left(2^{1/3}s\right)G_0(s)+2=0,
\end{equation}
where $F_0(s)$ is given by eqn.~(\ref{form:F0dyck}). On the level of 
coefficients, we have the recursion 
\begin{equation}\label{recG}
\alpha_{k-1}g_{k-1}+\sum_{l=0}^k2^{-l/2+5/2}f_lg_{k-l}=0,
\end{equation}
with the initial condition $g_0=s_{0,0,0}=2^{-1/2}$. If we set
\begin{displaymath}
g_k=\frac{\omega_k}{2^{3k/2+1/2}},
\end{displaymath}
then the above recursion is identical to that occurring in the definition 
of the meander distribution. In particular, all numbers $g_k$ are non-zero. 
Since the functions $S_{k,0}(x)$ are $\Delta(1/2)$-regular, we may use the 
transfer lemma \cite[Thm.~1]{FO90} to infer for the moments of $X_m^{(sym)}$ 
the asymptotic form
\begin{displaymath}
\begin{split}
\frac{\mathbb E_m[(X_m^{(sym)})_k]}{k!} &= 
\frac{1}{k!}\frac{[x^m]S_{k,0}(x)}{[x^m]S_{0,0}(x)}
\sim \frac{1}{k!}\frac{[x^m]s_{k,0,0}(1-2x)^{-(3k/2+1/2)}}
{[x^m]s_{0,0,0}(1-2x)^{-1/2}}\\ 
&\sim \frac{g_{k}}{g_{0}}
\frac{\Gamma(1/2)}{\Gamma(3k/2+1/2)} m^{3k/2}
= \frac{1}{2^{k}}\frac{\omega_k}{\omega_0}\frac{\Gamma(\alpha_0)}{\Gamma(\alpha_k)}
\frac{1}{2^{k/2}}m^{3k/2}
\qquad (m\to\infty).
\end{split}
\end{displaymath}
The last term is, up to the factor $m^{3k/2}$, the $k$-th moment of $Z/2$,
where $Z$ is the meander area variable. The previous estimate shows that 
the factorial moments $\mathbb E_m[(X_m^{(sym)})_k]$ are asymptotically 
equal to the ordinary moments $\mathbb E_m[(X_m^{(sym)})^k]$. It follows 
with \cite[Thm.~4.5.5]{Chu74} that the sequence of random variables $\left 
\{2m^{-3/2}X_m^{(sym)}\right \}_{m\in \mathbb N}$ converges in distribution 
to $Z$, where $Z$ is distributed as the meander area. We also have moment 
convergence.

\end{proof}

\begin{remark} \rm
As for the full class of staircase polygons, corrections to the asymptotic 
behaviour of the factorial moment generating functions can be mechanically 
obtained also for this example, by expanding the corresonding functional 
equation to higher orders in $\epsilon$.
\end{remark}

\subsection{Limit law for $\langle d_1\rangle$-symmetric polygons}\label{limitD1}

These polygons always have even half-perimeter. In order to derive an area limit 
law, we thus restrict to area random variables of quarter-perimeter. 
Note that we have $\widetilde D_1(x,q)=D_1(x^{1/2},q)$ for 
the generating function of the class $\langle d_1\rangle$-symmetric polygons, counted 
by quarter-perimeter and area. The functional equation for $D_1(x,q)$ induces
a similar one for $\widetilde D_1(x,q)$. Their factorial moment generating functions 
all have radius of convergence $1/4$, and a statement as in Lemma \ref{expDkl} 
can be formulated and proven almost verbatim for $\widetilde D_1(x,q)$. The method 
of dominant balance then yields a generating function for the leading coefficients 
in the singular expansions, a defining equation similar to eqn.~\eqref{diffeqF}, 
and a recursion similar to eqn.~\eqref{recF}. We have the following result.
\begin{theorem}
The area random variables $X_m$ of $\langle d_1\rangle$-symmetric 
staircase polygons, indexed by quarter-perimeter $m$ and scaled by 
$m^{-3/2}$, converge in distribution to a random variable $Y$, which 
is Airy distributed. We also have moment convergence. 
\end{theorem}
\subsection{Limit law for $\langle d_2\rangle$-symmetric polygons}

In \cite{LR01}, a combinatorial bijection between $\langle d_2 \rangle$-symmetric 
polygons and $\langle r^2\rangle$-symmetric polygons with even half-perimeter 
is described: cut a $\langle d_2 \rangle$-symmetric polygon along the line of 
reflection, flip its upper right part, and glue the two parts together along the 
cut. So Theorem \ref{theo:r2} translates to the $\langle d_2\rangle$-case.

Alternatively, one may apply the methods of Section~\ref{limitS}, together with 
modifications similar to those of Section~\ref{limitD1}, to the quarter-perimeter and
area generating function $\widetilde D_2(x,q)=D_2(x^{1/2},q)$.  
Lemma \ref{lem:symdyck} holds in this case, with $1/2$ replaced by $1/4$, and the
method of dominant balance yields results similar to eqn.~\eqref{diffeqG} and 
eqn.~\eqref{recG}.

\begin{theorem}
The area random variables $X_m$ of $\langle d_2\rangle$-symmetric 
staircase polygons, indexed by quarter-perimeter $m$ and scaled by 
$(2m)^{-3/2}$, converge in distribution to a random variable $Z/2$, 
where $Z$ is the meander area random variable. We also have moment 
convergence.
\end{theorem}

\subsection{Limit law for $\langle d_1,d_2 \rangle$-symmetric polygons}
In this symmetry class, every polygon has even half-perimeter. 
So we define $\widetilde D_{12}(x,q)=D_{12}(x^{1/2},q)$ as above, and obtain from 
the functional equation for $D_{12}(x,q)$ one for $\widetilde D_{12}(x,q)$,
which involves $\widetilde D_1(x,q)$, resembling eqn.~\eqref{feqS}. 

It can be argued, as in the proof Lemma \ref{lem:symdyck}, that all factorial
moment generating functions 
$\left.\frac{\partial^{k}}{\partial q^k}\widetilde D_{12}(x,q)\right|_{q=1}$ have radius 
of convergence $1/2$, with singularities at  $\pm1/2$,
where the leading singular behaviour of the coefficients is determined by the
singularity at $1/2$. We can apply the methods of Section~\ref{limitS}, with 
the modifications of Section~\ref{limitD1}. This yields the following result.

\begin{theorem}
The sequence of area random variables $X_m$ of $\langle d_1,d_2 
\rangle$-symmetric staircase polygons, indexed by quarter-perimeter $m$ and 
scaled by $m^{-3/2}$, converges in distribution to $2Z$, where $Z$ is
the meander area variable. We also have moment convergence.
\end{theorem}

\subsection{Limit law for $\langle r \rangle$-symmetric polygons}

The class of staircase polygons with $\langle r \rangle$-symmetry is the 
class of squares. These may be counted by quarter-perimeter $m$. Since 
for given $m$ there is exactly one square, they have, after scaling by 
$m^{-2}$, a concentrated limit distribution $\delta(x-1)$. This result 
can also be obtained from the $q$-difference equation in 
Proposition~\ref{prop:feqns}.

\subsection{Limit law for $\langle h,v \rangle$- ($\langle h \rangle$-,
$\langle v \rangle$-) symmetric polygons}

The class of staircase polygons with $\langle h,v \rangle$-symmetry (or with 
$\langle h \rangle$- or $\langle v \rangle$-symmetry) is the class of rectangles. 
These have been discussed in \cite{R07}.
The $k$-th moments of the area random variable $X_m$, with $m$ half-perimeter,
cf. eqn.~\eqref{tXm}, are given explicitly by
\begin{equation*}
\begin{split}
\mathbb E [X_m^k] = \sum_{l=1}^{m-1} (l(m-l))^k \frac{1}{m-1}
\sim m^{2k} \int_0^1 (x(1-x))^k {\rm d}x
=\frac{(k!)^2}{(2k+1)!}m^{2k} \qquad (m\to\infty),
\end{split}
\end{equation*}
where we used a Riemann sum approximation.
Consider the normalised random variable 
\begin{equation*}
\widetilde X_m=4 X_m/m^2.
\end{equation*}
The moments of $\widetilde X_m$ converge as $m\to\infty$, and the 
limit sequence $M_k=\lim_{m_\to\infty} \mathbb E_m[\widetilde X_m^k]$
satisfies Carleman's condition and hence defines a unique random 
variable with moments $M_k$. The corresponding distribution is 
the beta distribution $\beta_{1,1/2}$. We arrive at the following 
result.

\begin{theorem}
The sequence $\widetilde X_m=4 X_m/m^2$ of area random variables of 
rectangles, with half-perimeter $m$ scaled by $4/m^2$, converges 
in distribution to a $\beta_{1,1/2}$-distributed random variable. 
We also have moment convergence. 
\end{theorem}

One may also obtain this result by manipulating the associated 
$q$-difference equation, see \cite{R07}. Expansions of the 
factorial moment generating functions about their singularity at 
$x=1$ can be derived, and bounds for their most singular exponent 
can be given. The method of dominant balance can then be applied 
to obtain the leading singular coefficient of these 
expansions.

\section{Limit law for orbit counts}

Let $\mathcal{H}$ be a subgroup of $\mathcal{D}_4$. By the
Lemma of Burnside, the half-perimeter and area generating function
$P_\mathcal{H}(x,q)$ of orbit counts w.r.t. $\mathcal{H}$ is given by 
\cite{LR01} 
\begin{displaymath}
P_{\mathcal{H}}(x,q)=\frac{1}{|\mathcal{H}|}\sum_{g\in\mathcal{H}} 
P_{\mathrm{Fix}(g)}(x,q)
= \frac{1}{|\mathcal{H}|}\left( P(x,q)+R(x,q)\right),
\end{displaymath}
where $|\mathcal{H}|$ denotes the cardinality of $\mathcal{H}$,
and where $\mathrm{Fix}(g)\subseteq \mathcal{P}$ is the subclass of 
staircase polygons, which are fixed under $g\in\mathcal{D}_4$, with 
half-perimeter and area generating function $P_{\mathrm{Fix}(g)}(x,q)$.
The series $P(x,q)$ is the full staircase polygon half-perimeter and
area generating function, and $R(x,q)$ is the sum of generating
functions of staircase polygons which are fixed under $g\in\mathcal{H}$, where $g\ne e$.
Let $P(x,1)=\sum_m p_m x^m$ and  $R(x,1)=\sum_m r_m x^m$. Due
to the previous discussion, see also \cite[Prop.~14]{LR01}, the 
number of polygons fixed by some non-trivial symmetry grows subexponentially 
w.r.t. the total number of polygons. This implies
\begin{equation}\label{eqn:subex}
\frac{m^\alpha r_m}{p_m}\to 0 \qquad (m\to\infty),
\end{equation}
for any real number $\alpha$. 
As a consequence, area limit distributions for orbit
counts coincide with those for the full class of staircase
polygons. 

\begin{theorem}
Let $\mathcal{H}$ be a subgroup of $\mathcal{D}_4$. Then the area limit law
of the class $\mathcal{P}/\mathcal{H}$ coincides with that of $\mathcal{P}$.
\end{theorem}

\begin{proof}
We show that both classes have asymptotically the same area moments. Since in all
examples the limit distribution is uniquely determined by its moments,
the claim follows.

\smallskip

Note that, for polygons of half-perimeter $m$, their area $n$ satisfies 
$1\le n\le m^2$. Let $P(x,q)=\sum p_{m,n}x^mq^n$ and 
$R(x,q)=\sum r_{m,n}x^mq^n$. By eqn.~\eqref{eqn:subex}, we 
have for $k\in\mathbb N_0$
\begin{displaymath}
\frac{\sum_n n^k r_{m,n}}{\sum_n n^k p_{m,n}}\le
\frac{m^{2k}r_m}{p_m}\to 0 \qquad (m\to\infty).
\end{displaymath}
This implies for the coefficients of the moment generating functions 
the asymptotic estimate
\begin{displaymath}
\begin{split}
[x^m]\left(q\frac{\partial}{\partial q} \right)^k& \left.
P_\mathcal{H}(x,q)\right|_{q=1}=\frac{1}{|\mathcal{H}|}\sum_n n^k (p_{m,n}+r_{m,n})\\
&\sim \frac{1}{|\mathcal{H}|}\sum_n n^k p_{m,n}=
\frac{1}{|\mathcal{H}|}[x^m]\left.\left(q\frac{\partial}{\partial q} \right)^k 
P(x,q)\right|_{q=1} \qquad (m\to\infty).
\end{split}
\end{displaymath}
We conclude that both classes have asymptotically the same area moments. 
\end{proof}

\section{Staircase polygons, Dyck paths, and Brownian excursions}
\label{sec:stoch}

We briefly explain how some of our results could be alternatively obtained
from the bijections described in \cite{LR01}, which set symmetry classes of staircase 
polygons in one-to-one correspondence to symmetry classes 
of Dyck paths. Random Dyck paths, in turn, are 
related to corresponding stochastic objects such as the Brownian excursion 
and the Brownian meander. This allows to infer area limit distributions for 
some polygon symmetry classes from distributions of certain Brownian 
excursion and meander functionals. 

\smallskip

A combinatorial bijection between staircase polygons of perimeter $2m+2$
and Dyck paths of length $2m$ has been described by Delest and
Viennot \cite{DV84}. Within that bijection, the area of a staircase 
polygon corresponds to the sum of the peak heights of a Dyck path. 
In a uniform ensemble,
the sequence of random Dyck paths w.r.t.~half-length yields, 
after suitable normalisation, a sequence of stochastic processes 
on $(C[0,1],||\cdot||_\infty)$ with the Borel $\sigma$-algebra, 
which converges in distribution to the standard Brownian 
excursion \cite{B72, A91}. Also,
for certain continuous functionals on suitably normalised Dyck paths, 
including the area and the
sum of peak heights, we have convergence in distribution
and moment convergence to the Brownian excursion area, see 
\cite{GP92,G99} and \cite[Thm.~9]{D04}.
Using stochastic techniques, the distribution of the excursion 
area has been initially analysed in \cite{CR73,S82,L84}, compare 
the historical remarks in \cite{PW95}. Derivations by 
discrete methods are summarised in \cite{FL01}.

The above bijection, restricted to $r^2$-symmetric polygons,
yields \emph{symmetric} Dyck paths, which decompose in two identical
discrete meanders \cite{LR01}. Again, it is known that the sequence of 
random discrete meanders w.r.t.~half-length, appropriately normalised,
converges to the Brownian meander \cite{I74}. Together with
the convergence theorem for continuous functionals of polynomial 
growth \cite[Thm.~9]{D04}, one may conclude that also 
discrete meander functionals, such as area and 
sum of peak heights, converge to the Brownian meander area. 
The distribution of the Brownian meander area was derived in 
\cite{T95} from a $q$-difference equation similar to eqn.~\eqref{feqS}. 
It can also be obtained from known results for the corresponding
Brownian motion and bridges distribution, since there is a relation 
between the double Laplace transforms of the three distributions 
\cite{RY98, J07}. The same considerations hold for $d_2$-symmetric 
polygons, which are in one-to-one correspondence to $r^2$-symmetric 
polygons \cite{LR01}.

Staircase polygons with $d_1$-symmetry are in bijection
with pairs of identical Dyck paths, which is seen by cutting
a polygon along its positive diagonal \cite{LR01}. Here, the polygon 
area corresponds to twice the Dyck path area, hence
the area limit distribution is given by that of the Brownian
excursion area, by the convergence result mentioned above. 
A similar result holds for $\langle d_1, d_2\rangle$-symmetric 
polygons. Here, a combinatorial bijection between polygons and discrete 
meanders is known \cite{LR01}, where the polygon area corresponds 
to four times the meander area.

The classes of rectangles and squares lead to Dyck paths
with an initial sequence of up steps, followed by an alternating
sequence of peaks and valleys, ending in a terminal sequence of
down steps. As seen above, these classes are treated by
elementary methods.

\section{Conclusions}

We analysed the symmetry subclasses of 
staircase polygons on the square lattice. Exploiting
a simple decomposition for staircase polygons \cite{R06},
we obtained the area limit laws in the uniform
fixed perimeter ensembles.
This extends and completes previous results \cite{LR01}. 
As expected, orbit counts with respect to different 
symmetry subgroups always lead to an Airy distribution. 
The enumeration of polygons fixed under 
a given symmetry group leads to a 
variety of area limit distributions, such as a concentrated 
distribution, the $\beta_{1,1/2}$-distribution, 
the Airy-distribution, or the Brownian meander area distribution.

As described in Section~\ref{sec:stoch}, the latter two results
can also be obtained from the connection to Brownian 
motion and the Brownian meander. Our independent discrete approach 
uses only elementary methods from probability and singularity analysis of 
generating functions. Moreover, it may be used in order to analyse 
corrections to the limiting behaviour, which cannot easily be obtained
by stochastic methods.

One may also study the analogous problem of
\emph{perimeter} limit laws in a uniform ensemble where,
for fixed area, every polygon occurs with the same
probability. For the class of staircase polygons,
the associated centred and normalised random variable
is asymptotically Gaussian \cite[Prop.~9.11]{FS07}, and the same
result is expected to hold for the symmetry subclasses,
apart from squares. Also, limit laws in non-uniform ensembles and
for other counting parameters may be studied, compare \cite{R07}.

The above methods may be applied to extract limit laws
related to symmetry subclasses of other polygon 
classes. In particular, the classes of convex polygons 
on the square \cite{LRR98} and on the hexagonal \cite{GL05} 
lattices may be studied.

With respect to symmetry subclasses of self-avoiding polygons, 
exact enumeration studies may be carried out. It would be interesting
to numerically analyse moments, in order to test conjectures for 
limit distributions, compare  \cite{RGJ01,R07}.

\section*{Acknowledgements}

US and BT would like to acknowledge financial support by the German Research
Council (DFG) within the CRC701.

\end{document}